\documentclass[12pt]{article}
\usepackage{amssymb,amsmath,amscd}
\usepackage{theorem}

\newcommand {\book}[5]{{\sc #1:} {``#2."} {#3,} {#4} {(#5).}}

\newcommand \subs {\subseteq}
\newcommand \cO {{\mathcal {O}}}
\newcommand \bbF {{\mathbb {F}}}
\newcommand \eps {\varepsilon}

\theorembodyfont{\rmfamily}

\newcommand \proof {{\sc{Proof}}:\hspace*{3mm}}

\def\eop{{\unskip\nobreak\hfil\penalty50\hskip8mm\hbox{}
  \nobreak\hfil
  {$\dashv$}\parfillskip=0mm \par\smallskip}}

\begin{document}

\begin{center}
      {\Large{\sc{On The Complexity Of Hamel Bases\\[1.8ex]
                  Of Infinite Dimensional Banach Spaces}}}
\end{center}
\smallskip

\begin{center}
{\sc{Lorenz Halbeisen\footnote{I would like to thank the {\it Swiss
National Science Foundation\/} for its support during the period in
which the research for this paper has been done.}}}\\[1.7ex]
{\small{\sl Department of Pure Mathematics\\ Queen's University
Belfast\\ Belfast BT7 1NN, Northern Ireland}}\\[1ex] {\small{\sl
Email: halbeis@qub.ac.uk}}
\end{center}
\medskip

\hfill\smallskip

{\small {\bf 2000 Mathematics Subject Classification:} {\bf 46B20},
54E52.} \hfill\bigskip


\begin{abstract}
We call a subset $S$ of a topological vector space $V$ {\it linearly
Borel}, if for every finite number $n$, the set of all linear
combinations of $S$ of length $n$ is a Borel subset of $V$. It will
be shown that a Hamel base of an infinite dimensional Banach space
can never be linearly Borel. This answers a question of Anatolij
Plichko.
\end{abstract}

In the sequel, let $X$ be any infinite dimensional Banach space. A
subset $S$ of $X$ is called {\bf linearly Borel (w.r.t.\;$X$)}, if
for every positive integer $n$, the set of all linear combinations
with $n$ vectors of $S$ is a Borel subset of $X$. Since $X$ is a
complete metric space, $X$ is a {\bf Baire space}, {\sl i.e.}, a
space in which non-empty open sets are not meager ({\sl
cf.}\;\cite[Section\,3.9]{hitchhiker}). Moreover, all Borel subsets
of $X$ have the {\bf Baire property}, {\sl i.e.}, for each Borel set
$S$, there is an open set $\cO$ such that $\cO\Delta S$ is meager,
where $\cO\Delta S =(\cO\setminus S)\cup (S\setminus\cO)$.

This is already enough to prove the following.

{\sc Theorem.} If $X$ is an infinite dimension Banach space and $H$
is a Hamel base of $X$, then $H$ is not linearly Borel (w.r.t.\;$X$).

\proof Let $X$ be any infinite dimensional Banach space over the
field $\bbF$ and let $H$ be any Hamel base of $X$. For a positive
integer $n$, let $[H]^n$ be the set of all $n$-element subsets of $H$
and let $$H_n:=\bigg{\{} \sum_{i=1}^n \alpha_i h_i:\;
\alpha_1,\ldots,\alpha_n\in\bbF\setminus\{0\}\ \text{\sl{and}\/}\
\{h_1,\ldots,h_n\} \in [H]^n\bigg{\}}\,.$$ Assume towards a
contradiction that $H$ is linearly Borel. Then, by definition, for
each positive integer $n$, $H_n$ is Borel, and hence, by the facts
mentioned above, $H_n$ has the Baire property. Since $H$ is a Hamel
base of $X$, we get $$B=\bigcup_{n=1}^{\infty}H_n\,,$$ and because
$X$ is a Baire space, there must be a least positive integer $m$ such
that $H_m$ is not meager. Because $H_m$ has the Baire property and is
not meager, there is a non-empty open set $\cO$ such that $\cO\Delta
H_m$ is meager. Since $H$ is a Hamel base, $\cO\setminus H_m$ cannot
be empty, and therefore, $\cO\setminus H_m$ is a non-empty meager
set. Let $B_{v,r}$ denote the open ball with center $v\in X$ and
radius $r$. Let $x\in H_m\cap\cO$ and let $\eps$ be such that
$B_{x,2\eps}\subs\cO$. Take any $y\in H_{3m+1}$ with $\| y\|<\eps$,
then $B_{x+y,\eps}\subs\cO$. The following map is a homeomorphism
from $B_{x,2\eps}$ into $B_{x+y,\eps}$:
 $$\begin{array}{rrcl}
   \Phi : & B_{x,2\eps} & \longrightarrow & B_{x+y,\eps} \\
          & z           & \longmapsto     & x+y+\frac{1}{2}(z-x)
   \end{array}$$
Since $\cO\setminus H_m$ is meager, both sets, $B_{x,2\eps}\setminus
H_m$ as well as $B_{x+y,\eps}\setminus H_m$, are meager, and further,
by the definition of $\Phi$, also $B_{x+y,\eps}\setminus\Phi[H_m]$ is
meager, where $\Phi[H_m]:= \{\Phi(z):z\in H_m\cap B_{x,2\eps}\}$.
Now, because we have chosen $y\in H_{3m+1}$, $\Phi[H_m]\cap H_m
=\emptyset$, and hence, $$B_{x+y,\eps}=\big(B_{x+y,\eps}\setminus
H_m\big)\cup \big(B_{x+y,\eps}\setminus \Phi[H_m]\big)\,,$$ which
implies that the open set $B_{x+y,\eps}$, as the union of two meager
sets, is meager. But this is a contradiction to the fact that $X$ is
a Baire space. \eop

\end{document}